\documentclass[a4paper,11pt,english]{article}
\usepackage{amsmath}
\usepackage{amsfonts}
\usepackage{amssymb}
\usepackage{amsthm}
\usepackage[all]{xy}
\usepackage{geometry}
\usepackage{babel}

\newtheorem {thm}{Theorem}
\newtheorem* {thm*}{Theorem}
\newtheorem*{main*}{Main Theorem}

\newtheorem {cor}[thm]{Corollary}
\newtheorem {lem}[thm]{Lemma}
\newtheorem {prop}[thm]{Proposition}

\theoremstyle{definition}
\newtheorem {rem}[thm]{Remark}

\newtheorem {exa}[thm]{Example}

\DeclareMathOperator{\ord}{ord}

\DeclareMathOperator{\End}{End}

\newcommand{\hatj}{\hat{\jmath}}
\DeclareMathOperator{\Gal}{Gal}
\DeclareMathOperator{\Hom}{Hom}
\DeclareMathOperator{\lcm}{l.c.m.}

\title{Two variants of the support problem for products of abelian varieties and tori}
\author{Antonella Perucca}
\date{}
\begin{document}
\maketitle

\begin{abstract}
Let $G$ be the product of an abelian variety and a torus defined over a number field $K$. Let $P$ and $Q$ be $K$-rational points on $G$. 
Suppose that for all but finitely many primes $\mathfrak p$ of $K$ the order of $(Q \bmod \mathfrak{p})$ divides the order of $(P \bmod \mathfrak{p})$.
Then there exist a $K$-endomorphism $\phi$ of $G$ and a non-zero integer $c$ such that $\phi(P)=cQ$.
Furthermore, we are able to prove the above result with weaker assumptions: instead of comparing the order of the points we only compare the radical of the order (radical support problem) or the $\ell$-adic valuation of the order for some fixed rational prime $\ell$ ($\ell$-adic support problem). 
\end{abstract}

\section{Introduction}

Let $G$ be the product of an abelian variety and a torus defined over a number field $K$. Let $R$ be a $K$-rational point on $G$ and let $\phi$ be a $K$-endomorphism of $G$.
Then for all but finitely many primes $\mathfrak p$ of $K$ the order of $(\phi(R) \bmod \mathfrak{p})$ divides the order of $(R \bmod \mathfrak{p})$.
The support problem is concerned with the converse: what can we say about two $K$-rational points $P$ and $Q$ satisfying the following condition?

\begin{description}
\item[(SP)] The order of $(Q \bmod \mathfrak{p})$ divides the order of $(P \bmod \mathfrak{p})$ for all but finitely many primes $\mathfrak p$ of $K$.
\end{description}

This question was first studied in \cite{CorralesSchoof}, \cite{Kharepreprint} and \cite{BGKjacobians}.
Larsen solved the support problem for abelian varieties by showing that there exist a $K$-endomorphism $\phi$ and a non-zero integer $c$ such that $\phi(P)=cQ$ (\cite[Theorem 1]{Larsen03}). In general, one can not take $c=1$ even if $P$ and $Q$ have infinite order (\cite[Proposition~2]{Larsen03}).

We study two variants of the support problem, which we call respectively \emph{$\ell$-adic support problem} and \emph{radical support problem}. We require weaker conditions on the points:

\begin{description}
\item[(LSP)] Fix a rational prime $\ell$ and suppose that the $\ell$-adic valuation of the order of $(Q \bmod \mathfrak{p})$ is less than or equal to the $\ell$-adic valuation of the order of $(P \bmod \mathfrak{p})$, for all but finitely many primes $\mathfrak p$ of $K$.

\item[(RSP)] Fix an infinite set $S$ of rational primes and suppose that for every $\ell$ in $S$ the order of $(Q \bmod \mathfrak{p})$ is coprime to $\ell$ whenever the order of $(P \bmod \mathfrak{p})$ is coprime to $\ell$, for all but finitely many primes $\mathfrak p$ of $K$.
\end{description}

We strengthen Larsen's result on the support problem by proving the following:

\begin{main*}\label{maintwo}
Let $G$ be the product of an abelian variety and a torus defined over a number field $K$. Let $P$ and $Q$ be $K$-rational points on $G$. 
Suppose that $P$ and $Q$ satisfy condition (LSP) or condition (RSP).
Then there exist a $K$-endomorphism $\phi$ of $G$ and a non-zero integer $c$ such that $\phi(P)=cQ$.
\end{main*}

For abelian varieties, our result has an alternative proof: the proof by Larsen of \cite[Theorem 1]{Larsen03} only requires condition (RSP); the proof by Wittenberg of \cite[Theorem 1]{Larsen03} inspired from \cite{Larsenwhitehead} only requires condition (LSP), see \cite{Wittenberg}. For the multiplicative group or simple abelian varieties and assuming condition (LSP), equivalent results were proven respectively by Khare in \cite[Proposition 3]{Kharegalois} and by Bara\'nczuk in \cite[Theorem 8.2]{Baranczuk06}.
\newline

Let $G$ be the product of an abelian variety and a torus defined over a number field $K$. Let $P$ and $Q$ be points in $G(K)$ satisfying one of the conditions above. Let $c$ be the least positive integer such that $cQ$ belongs to the left $\End_K G$-submodule of $G(K)$ generated by $P$.
We prove the following:

Assuming condition (SP), $c$ divides a non-zero integer $m$ which depends only on $G$ and $K$. For abelian varieties this result has an alternative proof by Larsen, see \cite{Larsenwhitehead}.

Assuming condition (LSP), the $\ell$-adic valuation of $c$ is less than or equal to the $\ell$-adic valuation of a non-zero integer $m$ which depends only on $G$ and $K$ (notice that $m$ does not depend on $\ell$).

Assuming condition (RSP), there exists a non-zero integer $m$ depending only on $G$ and $K$ such that the following holds: for every $\ell$ in $S$ coprime to $m$ the $\ell$-adic valuation of $c$ is zero.

See section~\ref{sectionallc} for more results concerning $c$ under conditions (SP), (LSP) and (RSP) respectively.
\newline

Finally we discuss the \emph{multilinear support problem}, which is a variant of the support problem introduced by Bara\'nczuk in \cite{Baranczuk06}. The points $P$ and $Q$ are replaced by $n$-tuples of points and the following condition is required:

\begin{description}
\item[(MSP)] Suppose that for all but finitely many primes $\mathfrak p$ of $K$ and for all positive integers $m_1,\ldots,m_n$ the point $(m_1Q_1+\ldots+m_nQ_n \bmod \mathfrak p)$ is zero whenever the point $(m_1P_1+\ldots+m_nP_n \bmod \mathfrak p)$ is zero.
\end{description}

This condition is stronger than requiring condition (SP) on each pair of points $(P_i, Q_i)$ so there exist $K$-endomorphisms $\phi_i$ and an integer $c$ such that $\phi_i(P_i)=cQ_i$. One would like to prove that $\phi_i$ and $\phi_j$ are related for $i\neq j$. This is true if the endomorphism ring is $\mathbb Z$ (see \cite{Baranczuk06}) but in general $\phi_i$ and $\phi_j$ are not related for $i\neq j$, see section~\ref{sectionmultilinear}.
Another multilinear condition has recently been considered by Bara\'nczuk, see \cite{Baranczuk08}.

\section{Preliminaries}

Let $G$ be the product of an abelian variety and a torus defined over a number field $K$. Let $R$ be a $K$-rational point on~$G$ and call $G_R$ the smallest algebraic $K$-subgroup of $G$ containing $R$. 
Write $G_R^0$ for the connected component of the identity of $G_R$ and write $n_R$ for the number of connected components of $G_R$.
 
We say that $R$ is \textit{independent} if $R$ is non-zero and $G_R=G$.
The point $R$ is independent in $G$ if and only if $R$ is independent in $G\times_K \bar{K}$. Furthermore, $R$ is independent in $G$ if and only if $R$ is non-zero and the left $\End_K G$-submodule of $G(K)$ generated by $R$ is free. See \cite[Section 2]{Peruccaorder}.

\begin{prop}\label{ccbounded}
Let $G$ be the product of an abelian variety and a torus defined over a number field $K$.
Let $R$ be a $K$-rational point on $G$.
Then $n_R$ divides a non-zero integer which depends only on $G$ and $K$.
\end{prop}
\textit{Proof.} 
Write $G=A\times T$ and $R=(R_A,R_T)$. Since $G_R\subseteq G_{R_A}\times G_{R_T}$, we know that 
$G_R^0$ is the product of an abelian subvariety $A'$ of $G_{R_A}^0$ and a subtorus $T'$ of $G_{R_T}^0$ (see \cite[Proposition 5]{Peruccaorder}).
We have $A'=G_{R_A}^0$ because $A'$ contains a non-zero multiple of $R_A$. Analogously we have $T'=G_{R_T}^0$.
So  $G_R^0=G_{R_A}^0\times G_{R_T}^0$ hence $n_R$ divides the number of connected components of $G_{R_A}\times G_{R_T}$.
Then it suffices to prove the statement for $A$ and for $T$ respectively.

For $A$ the statement is proven in \cite[Lemma 2.2.4]{McQuillan}.
Now we prove the statement for $T$: we reduce at once to the case $T=\mathbb G_m^n$. 
Write $R=(R_1,\ldots,R_n)$ and let $e$ be the exponent of $\mathbb G_m(K)_{tors}$. Since $n_R$ divides $e$ times $n_{eR}$, we reduce to the case where $R_1,\ldots,R_n$ generate a torsion-free subgroup of $\mathbb G_m(K)$.
We conclude by proving that in this case $n_R=1$. We may clearly assume that $R$ is non-zero.
Fix a rational prime $\ell$.
Remark that $R_1,\ldots,R_n$ generate a free subgroup of $\mathbb G_m(K)$. By choosing a basis for this subgroup, we find an integer $s\geq 1$ and a point $R'$ independent in $\mathbb G_m^s$ such that  $\ord(R \bmod\mathfrak p)= \ord(R' \bmod\mathfrak p)$ for all but finitely many primes $\mathfrak{p}$ of $K$.
By \cite[Proposition 12]{Peruccaorder} there exist infinitely many primes $\mathfrak p$ such that $v_\ell[\ord(R' \bmod\mathfrak p)]=0$.
Then for infinitely many primes $\mathfrak p$ we have $v_\ell[\ord(R \bmod\mathfrak p)]=0$.  By \cite[Main Theorem]{Peruccaorder}, it follows that $v_\ell(n_R)=0$.
\hfill $\square$

\begin{lem}\label{extensionend}
Let $G$ be the product of an abelian variety and a torus defined over a number field $K$. Let $L$ be a finite Galois extension of $K$ of degree $d$.
Let $P$ and $Q$ be $K$-rational points on $G$. 
If $Q$ belongs to $\End_L G\cdot P$ then $dQ$ belongs to $\End_K G\cdot P$. \end{lem}
\textit{Proof.} 
Suppose that there exists $\psi$ in $\End_L G$ such that $\psi(P)=Q$.
Set $\phi=\sum_{\sigma\in\Gal(L/K)}\psi^\sigma$. 
Then $\phi$ is in $\End_K G$ and we have: $$\phi(P)=\sum_{\sigma\in\Gal(L/K)}\psi^\sigma(P)=\sum_{\sigma\in\Gal(L/K)}\psi(P)^\sigma=\sum_{\sigma\in\Gal(L/K)}Q^\sigma=dQ.$$
\hfill $\square$

\begin{lem}\label{isogeny} 
Let $A$ and $B$ be products of an abelian variety and a torus defined over a number field $K$.
Let $\alpha$ be an isogeny in $\Hom_K( A,B)$ and let $d$ be the exponent of the kernel of $\alpha$ (which divides the degree of $\alpha$).
Let $R$ be a $K$-rational point on $A$.
For all but finitely many primes $\mathfrak p$ of $K$ the following holds: the order of $(dR \bmod{\mathfrak p})$ divides the order of $(\alpha(R) \bmod{\mathfrak p})$. 
\end{lem}
\textit{Proof.} 
For every $\psi$ in $\Hom_K (B,A)$ and for every point $W$ in $B(K)$ the following holds: the order of $(\psi(W) \bmod{\mathfrak p})$ divides the order of $(W \bmod{\mathfrak p})$ for all but finitely many primes $\mathfrak p$ of $K$.
Call $\hat{\alpha}$ the isogeny in $\Hom_K (B,A)$ such that $\hat{\alpha}\circ\alpha=[d]$.
The statement follows by applying the first assertion to $\psi=\hat{\alpha}$ and $W=\alpha(R)$.
\hfill $\square$

\begin{lem}\label{lemnew}
Let $K$ be a number field. Let $I=\{1,\ldots,n\}$. For every $i\in I$ let $B_i$ be the product of an abelian variety and a torus defined over $K$. Suppose that for $i\neq j$ either $B_i=B_j$ or $\Hom_K(B_i,B_j)=\{0\}$. Let $H=\prod_{j\in J} B_j$ for some subset $J$ of $I$ and let $R$ be a point in $H(K)$ which is independent in $H$. Let $W$ be a point in $B_{n}(K)$. Then if $(R,W)$ is not independent in $H\times B_{n}$ there exists a non-zero $f$ in $\End_K B_{n}$ such that $f(W)$ belongs to $\Hom_K(H,B_{n})\cdot R$.
\end{lem}
\textit{Proof.} We know that there exists a non-zero $F$ in $\End_K(H\times B_{n})$ such that $F(R,W)=0$. 
Write $F=(F_1,F_2)$ in the decomposition $$\End_K(H\times B_{n})=\Hom_K(H, H\times B_{n})\times \Hom_K(B_{n}, H\times B_{n})\,. $$
We then have $F_1(R)+F_2(W)=F(R,W)=0$.

Since $F\neq 0$ there exists a factor $B_m$ of $H\times B_{n}$ such that $\pi_m\circ F\neq 0$ where $\pi_m$ is the projection of $H\times B_{n}$ onto $B_m$.
Now we prove that $\pi_m\circ F_2\neq 0$. Suppose not. Then we must have $\pi_m\circ F_1\neq 0$.
If $B_n$ is not equal to any factor of $H$ and $B_m=B_n$ we have $\Hom_K(H, B_{m})=\{0\}$ hence $\pi_m\circ F_1= 0$, contradiction.
So we may assume that there is an inclusion map $i$ from $B_m$ to $H$. We have $i\circ \pi_m\circ F_1\neq 0$ and $i\circ \pi_m\circ F_1(R)=-i\circ\pi_m\circ F_2(W)=0$, which contradicts the fact that $R$ is independent in $H$.

Since $\pi_m\circ F_2\neq 0$, we have $\Hom_K(B_{n},B_m)\neq \{0\}$ hence $B_m=B_{n}$. Call $f=\pi_m\circ F_2$. Then $f$ is a non-zero element of $\End_K B_{n}$ and we have  $f(W)=-\pi_m\circ F_1(R)$ hence $f(W)$ belongs to $\Hom_K(H,B_{n})\cdot R$.
 \hfill $\square$

\begin{lem}\label{estimate}
Let $K$ be a number field and let $I=\{1,\ldots,n\}$. Let $G=\prod_{i\in I} B_i$ where for every $i$ $B_i$ is either $\mathbb G_m$ or a $K$-simple abelian variety and for $i\neq j$ either $B_i=B_j$ or $\Hom_K(B_i,B_j)=\{0\}$.
Let $P=(P_1,\ldots,P_n)$ be a point on $G(K)$ of infinite order.
Then there exist a subset $J=\{j_1,\ldots, j_s\}$ of $I$ and a non-zero integer $d$ such that the point $P'=(P_{j_1},\ldots,P_{j_s})$ is independent in $G'=\prod_{j\in J} B_j$ and such that for all but finitely many primes $\mathfrak p$ of $K$ the order of $(P \bmod \mathfrak{p})$ divides $d$ times the order of $(P' \bmod \mathfrak{p})$.
\end{lem}
\textit{Proof.}
We prove the statement by induction on $n$. If $n=1$, the point $P_1$ is independent in $B_1$ so take $J=\{1\}$, $d=1$.
Now we prove the inductive step. Let $P=(P_1,\ldots, P_n)$ and set $\tilde{P}=(P_1,\ldots,P_{n-1})$. If $\tilde{P}$ is a torsion point then $P_n$ is independent in $B_n$ and we easily conclude. So assume that $\tilde{P}$ has infinite order and let $\tilde{J}$, $\tilde{d}$, $\tilde{P'}$ and $\tilde{G'}$ be as in the statement.
If the point $(\tilde{P'},P_{n})$ is independent in $\tilde{G'}\times B_n$ take $J=\tilde{J}\cup\{n\}$ and $d=\tilde{d}$.
Otherwise by Lemma~\ref{lemnew} there exists a non-zero $f$ in $\End_K B_{n}$ such that $f(P_n)$ is in $\Hom_K(\tilde{G'}, B_n)\cdot \tilde{P'}$.

Since $f$ is an isogeny, there exist $\hat{f}$ in $\End_K B_{n}$ and a non-zero integer $r$ such that $[r]=\hat{f}\circ f$.
Consequently $rP_n$ belongs to $\Hom_K(\tilde{G'}, B_n)\cdot \tilde{P'}$ and so we can take $J=\tilde{J}$ and $d=\lcm (\tilde{d},r)$.
\hfill $\square$

\begin{lem}\label{j}[Proposition 2, Appendix of \cite{Bertrandpreprint}]
Let $A$ be an abelian variety defined over a number field $K$. There exists a non-zero integer $t$ such that the following holds: for every $K$-rational point $R$ on $A$ there exists an abelian subvariety $Z$ of $A$ defined over $K$ such that $G^0_R+ Z=A$ and $G^0_R\cap Z$ has order dividing $t$.
\end{lem}

\noindent The previous lemma can also be found in \cite[ Proposition 5.1]{RatazziUllmo}.

\section{The proof of the Main Theorem}\label{sectionproof}

\begin{lem}\label{ladicisogenous}
Let $A$ and $B$ be products of an abelian variety and a torus defined over a number field $K$ and $K$-isogenous. If the Main Theorem is true for $B$, then it is true for $A$.
\end{lem}
\textit{Proof.} 
Suppose that the Main Theorem holds for $B$.
Let $\alpha$ be a $K$-isogeny from $A$ to $B$, call $d$ the degree of $\alpha$ and call $\hat{\alpha}$ the isogeny in $\Hom_K(B,A)$ satisfying $\hat{\alpha}\circ\alpha=[d]$.
Because of Lemma~\ref{isogeny},
if $P$ and $Q$ satisfy condition (LSP) then for all but finitely many primes $\mathfrak p$ of $K$ we have:
\begin{multline*}
v_\ell[\ord(\alpha(P)\bmod \mathfrak p)]\geq v_\ell[\ord(dP\bmod \mathfrak p)]
\geq v_\ell[\ord(dQ\bmod \mathfrak p)]\geq v_\ell[\ord(\alpha(dQ)\bmod \mathfrak p)].
\end{multline*}
So $\alpha(P)$ and $\alpha(dQ)$ satisfy condition (LSP).
By Lemma~\ref{isogeny}, if $P$ and $Q$ satisfy condition (RSP) then $\alpha(P)$ and $\alpha(Q)$ satisfy condition (RSP) for the subset of $S$ consisting of the primes coprime to $d$.
We deduce that $$\psi\bigl(\alpha(P)\bigr)=r\bigl(\alpha(dQ)\bigr)$$
where $\psi$ is in $\End_K B$ and $r$ is a non-zero integer.
Set $\phi=\hat{\alpha}\circ\psi\circ\alpha$, $c=rd^2$.
Then $\phi$ is in $\End_KA$, $c$ is a non-zero integer and we have: 
$$\phi\bigl(P\bigr)=\hat{\alpha}\circ \psi\circ \alpha\bigl(P\bigr)=\hat{\alpha}\circ [r]\circ \alpha \bigl(dQ\bigr)=rd^2Q=cQ.$$
\hfill $\square$\newline

\noindent \textit{Proof of the Main Theorem.} 

\noindent \emph{First step.}
We reduce to prove the theorem for $G=\prod_{i\in I} B_i$ where for every $i$ the factor $B_i$ is either $\mathbb G_m$ or a $K$-simple abelian variety and for $i\neq j$ either $B_i=B_j$ or $\Hom_K(B_i,B_j)=\{0\}$.
To accomplish this, it suffices to combine two things: the statement holds for $G$ if it holds for $G\times_K L$, where $L$ is a finite Galois extension of $K$; the statement holds for $G$ if it holds for $\alpha(G)$ where $\alpha$ is a $K$-isogeny.
The first assertion is a consequence of Lemma~\ref{extensionend}.
The second assertion is proven in Lemma~\ref{ladicisogenous}.\\
\noindent \emph{Second step.}
Let $G=\prod_{i\in I} B_i$ and write $P=(P_1,\ldots,P_n)$, $Q=(Q_1,\ldots,Q_n)$. Without loss of generality we may replace $Q$ by $(Q_1,0,\ldots,0)$.\\
We may assume that $Q$ has infinite order (otherwise take $\phi=0$ and $c=\ord Q$). Then we may assume that also $P$ has infinite order. Otherwise, let $\ell$ be either the prime of condition (LSP) or a prime of $S$ coprime to $\ord (P)$. We find a contradiction by \cite[Corollary 14]{Peruccaorder} since there exist infinitely many primes $\mathfrak p$ of $K$ such that $v_\ell[\ord(Q \bmod\mathfrak{p})]>v_\ell[\ord(P)]$.\\
\emph{Third step.}
Apply Lemma~\ref{estimate} to $P$ and let $J$, $d$, $P'$, $G'$ be as in Lemma~\ref{estimate}.
Since $P'$ is a projection of $P$, it suffices to prove that there exist $\psi$ in $\Hom_K(G',B_1)$ and a non-zero integer $c$ such that $\psi(P')=c Q_1$.\\
\emph{Fourth step.}
The point $(P',Q_1)$ is not independent in $G'\times B_1$.
Otherwise, let $\ell$ be either the prime of condition (LSP) or a prime of $S$ coprime to $d$ and apply \cite[Proposition 12]{Peruccaorder}.
There exist infinitely many primes $\mathfrak p$ of $K$ such that $v_\ell[\ord (P' \bmod \mathfrak p)]=0$ and $v_\ell[\ord (Q_1 \bmod \mathfrak p)]=v_\ell(d)+1$. We find a contradiction since by definition of $d$ we may assume that $v_\ell[\ord (P \bmod \mathfrak p)]\leq v_\ell(d)+ v_\ell[\ord (P' \bmod \mathfrak p)]$.\\
\emph{Fifth step.}
By definition $P'$ is independent in $G'$ so we can apply Lemma~\ref{lemnew} to the points $P'$ and $Q_1$. Then since $(P',Q_1)$ is not independent in $G'\times B_1$ there exists a non-zero $f$ in $\End_K B_1$ such that 
$f(Q_1)$ belongs to $\Hom_K(G',B_1)\cdot P'$.
Since $f$ is an isogeny, there exist $\hat{f}$ in $\End_K B_{1}$ and a non-zero integer $c$ such that $[c]=\hat{f}\circ f$.
Consequently $cQ_1$ belongs to $\Hom_K(G', B_1)\cdot P'$.
\hfill{$\square$}\newline

The following corollary is the analogue to \cite[Corollary 6]{Larsen03}.

\begin{cor}\label{thmduo}
Let $G_1$ and $G_2$ be products of an abelian variety and a torus defined over a number field $K$. Let $P$ and $Q$ be $K$-rational points respectively on $G_1$ and $G_2$ satisfying condition (LSP) or condition (RSP). Then there exist $\phi$ in $\Hom_K(G_1,G_2)$ and a non-zero integer $c$ such that $\phi(P)=cQ$.
\end{cor}
\textit{Proof.}
Apply the Main Theorem to $G_1\times G_2$ and the points $(P,0)$ and $(0,Q)$.
\hfill $\square$

\section{On the integer $c$ of the Main Theorem}\label{sectionallc}

The following proposition is the generalization of a result by Khare and Prasad (\cite[Theorem 1]{KharePrasad}).

\begin{prop}\label{constantindependent}
Under the assumptions of Corollary~\ref{thmduo} and if the point $P$ is independent in $G_1$, one can take $c$ coprime to $\ell$ under condition (LSP) and coprime to every $\ell$ in $S$ under condition (RSP).
\end{prop}
\textit{Proof.}
We have $\phi P=cQ$ for some $\phi$ in $\Hom_K(G_1,G_2)$ and some non-zero integer $c$.
Let $\ell$ be either the prime of condition (LSP) or a fixed prime of $S$.
By iteration, it suffices to prove that if $c$ is divisible by $\ell$ there exists $\psi$ in $\Hom_K(G_1,G_2)$ such that $\psi P=\frac{c}{\ell}Q$. 
So suppose that $c$ is divisible by $\ell$.
Let $P'$ be a point in $G_1(\bar{K})$ such that $\ell P'=P$.
We then have $\phi(P')=\frac{c}{\ell}Q+Z$ for some $Z$ in $G_2[\ell]$. 
Write $L$ for a finite extension of $K$ over which $G_1[\ell]$ is split and where $P'$ is defined.
Notice that $P'$ is also independent in $G_1$.
The condition of Corollary~\ref{thmduo} clearly implies that for all but finitely many primes $\mathfrak q$ of $L$ the order of $(Q \bmod\mathfrak q)$ is coprime to $\ell$ whenever the order of $(P \bmod\mathfrak q)$ is coprime to $\ell$.

First we prove that $\phi=[\ell]\circ\psi$ for some $\psi$ in $\Hom_K(G_1,G_2)$. Suppose not and then let $T$ be a point in $G_1[\ell]\backslash\ker(\phi)$. 

Suppose that $\phi(T)\neq Z$.
By \cite[Proposition 11]{Peruccaorder} there exist infinitely many primes $\mathfrak q$ of $L$ such that $v_\ell[\ord(P'-T \bmod \mathfrak q)]=0$. We deduce that $v_\ell[\ord(P \bmod \mathfrak q)]=0$ and that the point $(\phi(P')-\phi(T) \bmod \mathfrak q)$ has order coprime to $\ell$. Then 
$$r_{\mathfrak{q}}\phi(T)=r_{\mathfrak{q}}\phi(P')=r_{\mathfrak{q}}(\frac{c}{\ell}Q+Z) \;(\bmod \mathfrak{q})$$
for some integer $r_{\mathfrak{q}}$ coprime to $\ell$. Therefore 
$$r_{\mathfrak{q}}\frac{c}{\ell}Q=r_{\mathfrak{q}}(\phi(T)-Z)\;(\bmod \mathfrak{q}).$$
By discarding finitely many primes $\mathfrak q$, we may assume that the order of $(\phi(T)-Z \bmod \mathfrak q)$ is $\ell$. We deduce that $v_\ell[\ord(Q \bmod \mathfrak q)]>0$ and we find a contradiction.

Now suppose that $\phi(T)=Z$. Then $\phi(P')=\frac{c}{\ell}Q+\phi(T)$.
By \cite[Proposition 11]{Peruccaorder} there exist infinitely many primes $\mathfrak q$ of $L$ such that $v_\ell[\ord(P'\bmod \mathfrak q)]=0$. Then $v_\ell[\ord(P\bmod \mathfrak q)]=0$. 
By discarding finitely many primes $\mathfrak q$, we may assume that the order of $(\phi(T)\bmod \mathfrak q)$ is $\ell$. We deduce that 
$v_\ell[\ord(Q \bmod \mathfrak q)]>0$ and we find a contradiction.

So we can factor $\phi$ as $[\ell]\circ\psi$ for some $\psi$ in $\Hom_K(G_1,G_2)$.
Then $\psi(P)=\frac{c}{\ell}Q+T'$ for some $T'$ in $G_2[\ell]$. It suffices to prove that $T'=0$.
By \cite[Proposition 12]{Peruccaorder}, there exist infinitely many primes $\mathfrak q$ of $L$ such that $v_\ell[\ord(P \bmod \mathfrak q)]=0$. If $T'\neq 0$, we may assume that the order of $(T'\bmod \mathfrak q)$ is $\ell$. We deduce that $v_\ell[\ord(Q \bmod \mathfrak q)]>0$ and we have a contradiction.
\hfill $\square$\newline

\begin{prop}\label{boundedc}
Under the assumptions of the Main Theorem, let $c$ be the least positive integer such that $cQ$ belongs to $\End_K G\cdot P$.
If condition (LSP) holds then $v_\ell(c)\leq v_\ell(m)$ for some non-zero integer $m$ depending only on $G$ and $K$.
If condition (RSP) holds then $v_\ell(c)=0$ for every $\ell$ in $S$ coprime to $m$, for some non-zero integer $m$ depending only on $G$ and $K$.
\end{prop}
\textit{Proof.}
We first reduce to the case $G=A\times T$ where $A$ is an abelian
variety and $T=\mathbb G_m^n$.
It suffices to show that the statement holds for $G$ if it holds for $G\times_K L$ where $L$ is a finite Galois extension of $K$. This can be deduced from the proof of Lemma~\ref{extensionend}: if $m$ is as in the statement for $G\times_K L$ then for $G$ one can take $[L:K]m$.

We reduce to the case where $G_{P}$ is connected. By Proposition~\ref{ccbounded}, $n_P$ divides an integer $h$ depending only on $G$ and $K$. We can then replace $P$ and $Q$ with $hP$ and $hQ$.

If $P$ is zero then from \cite[Corollary 14]{Peruccaorder} we immediately deduce that $Q$ is a torsion point. In this case $c$ divides the exponent of $G(K)_{tors}$.

Now we assume that $G_{P}$ is connected and that $P$ has infinite order. 
By \cite[Proposition 5]{Peruccaorder}, we have $G_P=A'\times T'$ where $A'$ is an abelian subvariety of $A$ and $T'$ is a sub-torus of $\mathbb G_m^n$. 
Since $P$ is independent in $G_P$, from Proposition~\ref{constantindependent} it follows that there exist $\psi$ in $\Hom_K(G_{P},G)$  and an integer $r$ coprime to $\ell$ (respectively to every prime of $S$) such that $\psi(P)=r Q$.

Write $P=(P_A,P_T)$ and remark that $A'=G_{P_A}$ (see the proof of Proposition~\ref{ccbounded}).
Apply Lemma~\ref{j} to $P_A$. Let $Z$ and $t$ be as in Lemma~\ref{j}. Then the map $$j:A'\times Z\rightarrow A\,;\; (x,y)\mapsto x+y.$$ is a $K$-isogeny in $\Hom_K(A'\times Z,A)$ of degree dividing $t$.
Call $\hatj$ the isogeny in $\Hom_K(A, A'\times Z)$ satisfying $\hatj\circ j=[t]$. We have: $$\hatj\bigl(P_A\bigr)= \hatj\circ j\bigl((P_A,0)\bigr)=(tP_A,0).$$

Then there is an element $\pi_{A}$ in $\Hom_K(A,A')$ mapping $P_{A}$ to $tP_{A}$. Since $T'$ is a direct factor of $T$, there exists $\pi_T$ in $\Hom_K(T,T')$ such that $\pi_T(P_{T})=tP_{T}$. Let $\Pi$ be $\pi_A\times \pi_T$. Then $\Pi$ is in $\Hom_K(G,G_P)$ and $\Pi(P)=tP$.
The map $\phi=\psi\circ\Pi$ is in $\End_K G$ and we have $\phi(P)=rtQ$. 

Since $r$ is coprime to $\ell$ (respectively to every prime of $S$) and $t$ depends only on $G$ and $K$, this concludes the proof.
\hfill $\square$\newline

Unless $G(K)$ is finite, one clearly cannot bound $v_p(c)$ for any rational prime $p$ different from $\ell$ (assuming condition (LSP)) or not in $S$ (assuming condition (RSP)).

Assuming condition (LSP), a straightforward adaptation of \cite[Proposition 2]{Larsen03} shows that in general one cannot take $c$ coprime to $\ell$ even if $P$ and $Q$ have infinite order.

For a split torus or for an abelian variety and assuming condition (RSP), one cannot in general bound $v_\ell(c)$ for every $\ell$ in $S$:

\begin{exa}\label{radnoboundS}
Let $\ell$ be a rational prime. Let $G$ be either the multiplicative group or an elliptic curve without complex multiplication defined over a number field $K$. Suppose that $G(K)$ contains a point $R$ of infinite order and a torsion point $T$ of order $\ell$.
Consider the points $P=(\ell^hR,T)$ and $Q=(R,0)$ on $G^2$, for some fixed $h$ in $\mathbb N$. Then the points $P$ and $Q$ satisfy condition (RSP) where $S$ is the set of all primes but one has to take $c$ such that $v_\ell(c)\geq h$. 
By varying $h$, we see at once that that one cannot bound $v_\ell(c)$ with a constant depending only on $G$ and $K$.
\end{exa}

\begin{prop}\label{constanttori}
In the Main Theorem, assuming condition (LSP) and if $G$ is a split torus then one can take $c$ coprime to $\ell$. 
\end{prop}
\noindent \textit{Proof.}
We may assume that $G=\mathbb G_m^n$. Recall that $\mathbb G_m[a]\simeq \mathbb Z/a \mathbb Z$ for every $a\geq 1$.
Without loss of generality we may assume that $Q=(Q_1,0,\ldots,0)$.
If $P$ is a torsion point then (because of $\phi(P)=cQ$) $Q_1$ is also a torsion point and the statement easily follows from condition (LSP).
Now assume that $P$ has infinite order. Since $\End_K \mathbb G_m\simeq\mathbb Z$, we may assume that $P$ is of the following form:
$$P=(R_1,\ldots,R_h,T,0,\ldots,0)$$
where the point $(R_1,\ldots, R_h)$ is independent in $\mathbb G_m^h$, $h\geq 1$ and $T$ is a torsion point. Call $t$ the $\ell$-adic valuation of the order of $T$.

We have \begin{equation}\label{eqnmulti}
a T+\sum_{i=1}^h a_iR_i=cQ_1
\end{equation}
for some $a,a_1,\ldots,a_h$ in $\mathbb Z$ and for some non-zero integer $c$.
Suppose that $c$ is divisible by $\ell$. It suffices to find an expression analogous to \eqref{eqnmulti} where $c$ is replaced by $\frac{c}{\ell}$ and we conclude by iteration.

Now we prove that $a$ is divisible by $\ell$. Suppose not. We may clearly assume that $t\neq 0$, otherwise we can multiply every coefficient of \eqref{eqnmulti} by an integer coprime to $\ell$ and replace $a$ by zero.
By \cite[Proposition 12]{Peruccaorder} there exist infinitely many primes $\mathfrak p$ of $K$ such that $v_\ell[\ord(R_i \bmod \mathfrak p)]=0$ for every $i$. We may assume that $v_\ell[\ord(T \bmod \mathfrak p)]=t$.
We deduce that $v_\ell[\ord(Q \bmod \mathfrak p)]\geq t+1$ and that $v_\ell[\ord(P \bmod \mathfrak p)]= t$ so we find a contradiction.

Without loss of generality we prove that $a_h$ is divisible by $\ell$. Suppose not.
The point $(R_1,\ldots, R_{h-1}, a_hR_h+aT)$ is independent in $\mathbb G_m^h$. Thus by \cite[Proposition 12]{Peruccaorder} there exist infinitely many primes $\mathfrak p$ of $K$ such that $v_\ell[\ord(R_i \bmod \mathfrak p)]=0$ for every $i\neq h$ and $v_\ell[\ord(a_hR_h+aT \bmod \mathfrak p)]=t+1$.
We easily deduce that $v_\ell[\ord(Q \bmod \mathfrak p)]\geq t+2$ and that  $v_\ell[\ord(P \bmod \mathfrak p)]= t+1$, contradiction.

Now we can write $$\frac{a}{\ell} T+\sum_{i=1}^m \frac{a_i}{\ell}R_i=\frac{c}{\ell}Q_1+W$$
where $W$ is in $\mathbb G_m[\ell]$.

If $t\geq 1$ then $W$ is a multiple of $T$ and we conclude. If $W=0$ we also conclude.
Now suppose that $t=0$ and $W\neq 0$.
By \cite[Proposition 12]{Peruccaorder} there exist infinitely many primes $\mathfrak p$ of $K$ such that $v_\ell[\ord(R_i \bmod \mathfrak p)]=0$ for every $i$. We may assume that the order of $(W \bmod \mathfrak p)$ is $\ell$.
We deduce that $v_\ell[\ord(P \bmod \mathfrak p)]=0$ and  $v_\ell[\ord(Q \bmod \mathfrak p)]\geq 1$, a contradiction. 
\hfill $\square$\newline 

By the previous proposition and Lemma~\ref{extensionend}, assuming condition (LSP) for a  torus one can take $c$ such that $v_\ell(c)\leq v_\ell(d)$ where $d$ is the degree of a finite Galois extension of $K$ where the torus splits.
In particular, if $G$ is a $1$-dimensional torus one can take $c$ coprime to $\ell$ (since every endomorphism is defined over $K$).

We may weaken condition (LSP) in the Main Theorem as follows: there exists an integer $d\geq0$ such that for all but finitely many primes $\mathfrak p$ of $K$ $v_\ell[\ord(P\bmod \mathfrak p)]$ is greater than or equal to $v_\ell[\ord(Q \bmod \mathfrak{p})]-d$. Indeed, it is immediate to see that $P$ and $\ell^dQ$ satisfy condition (LSP).

Notice that the set $S$ in condition (RSP) needs in general to be infinite:

\begin{exa}
Let $S$ be a finite family of prime numbers and let $m$ be the product of the primes in $S$.
Let $G$ be either the multiplicative group or an elliptic curve without complex multiplication defined over a number field $K$. Suppose that $G(K)$ contains a torsion point $T$ of order $m$ and that the rank of $G(K)$ is greater than $1$. Then let $(R,W)$ be a point in $G^2(K)$ which is independent.
Consider the points $P=(R,T)$, $Q=(W,0)$ in $G^2(K)$.
The order of $P$ is a multiple of $m$ for all but finitely many primes $\mathfrak p$ of $K$ hence the points $P$ and $Q$ satisfy condition (RSP) for the set $S$. Nevertheless, no non-zero multiple of $Q$ lies in the left $\End_K G^2$-submodule of $G^2(K)$ generated by $P$.
\end{exa}

Now suppose that condition (SP) holds.
In general one can not take $c=1$ even if $P$ and $Q$ have infinite order (\cite[Proposition~2]{Larsen03}). As a consequence of Proposition~\ref{constantindependent}, one can take $c=1$ if $P$ is independent in $G$. This is the generalization of a result by Khare and Prasad (\cite[Theorem 1]{KharePrasad}).
As a consequence of Proposition~\ref{boundedc}, one can take $c$ such that it divides a constant depending only on $G$ and $K$. This was known for abelian varieties, see \cite[Corollary 4.4 and Theorem 5.2]{Larsenwhitehead} by Larsen. More precisely, Larsen proved that for abelian varieties one can take $c$ dividing the exponent of $G(K)_{tors}$ whenever the Tate-modules are all integrally semi-simple (and in every $K$-isogeny class there is such an abelian variety).
Notice that assuming condition (SP) it is not true in general that there exist a $K$-endomorphism $\phi$ and a $K$-rational torsion point $T$ such that $\phi(P)=Q+T$. A counterexample was found by Larsen and Schoof in  \cite{LarsenSchoof}.

\section{The multilinear support problem}\label{sectionmultilinear}

In this section we discuss the multilinear support problem, introduced by Bara\'nczuk in \cite{Baranczuk06}.
We first show that condition (MSP) (see the Introduction) is stronger than the condition of the support problem on each pair of points.

\begin{rem}\label{remmulti}
Assuming condition (MSP), the following holds: for every $i=1,\ldots, n$ the order of $(Q_i \bmod \mathfrak p)$ divides the order of 
$(P_i \bmod \mathfrak p)$ for all but finitely many primes $\mathfrak p$ of $K$.
\end{rem}
\begin{proof}
Without loss of generality it suffices to prove the claim for $P_1$ and $Q_1$. Let $\mathfrak p$ be a prime ideal of $K$ such that condition (MSP) holds. For every $i\neq1$ fix $m_i$ such that $(m_iP_i \bmod \mathfrak p)=0$ and $(m_iQ_i \bmod \mathfrak p)=0$. Then for every positive integer $m_1$ we have $(m_1Q_1 \bmod \mathfrak p)=0$ whenever $(m_1P_1 \bmod \mathfrak p)=0$. Consequently, the order of $(Q_1 \bmod \mathfrak p)$ divides the order of 
$(P_1 \bmod \mathfrak p)$.
\end{proof}

Because of the previous remark and the Main Theorem, there exist $K$-endomorphisms $\phi_i$ and an integer $c$ such that $\phi_i(P_i)=cQ_i$. One would like to prove that $\phi_i$ and $\phi_j$ are related for $i\neq j$. This is true if the endomorphism ring is $\mathbb Z$ (see \cite{Baranczuk06}) but in general $\phi_i$ and $\phi_j$ are not related for $i\neq j$:

\begin{exa}\label{notrelated}
Let $E$ be an elliptic curve defined over a number field $K$. Let $R_1$, $R_2$ be points in $E(K)$ and let $\phi_1$, $\phi_2$ be in $\End_K E$. The following points in $E^2(K)$ satisfy condition (MSP):
$$P_1=(R_1,0)\,;\;P_2=(0,R_2)\,;\; Q_1=(\phi_1(R_1),0)\,;\; Q_2=(0,\phi_2(R_2)).$$
\end{exa}

The next example shows that $\phi_i$ and $\phi_j$ are in general not related, not even for an elliptic curve, if we require the following weaker condition:

\begin{description}
\item[(LMSP)] Fix a rational prime $\ell$ and suppose that for all but finitely many primes $\mathfrak p$ of $K$ and for all positive integers $m_1,\ldots,m_n$ the order of  $(m_1Q_1+\ldots+m_nQ_n \bmod \mathfrak p)$ is coprime to $\ell$ whenever the order of $(m_1P_1+\ldots+m_nP_n \bmod \mathfrak p)$ is coprime to $\ell$.
\end{description}

\begin{exa}\label{ellipticcurveexa}
Let $E$ be an elliptic curve defined over a number field $K$ such that $\End_K E=\mathbb Z[\mathit i]$.
Let $\phi_1$ and $\phi_2$ be in $\End_K E$ and let $P_1$ be in $E(K)$. The following points satisfy condition (LMSP) for $\ell=3$:
$$P_1\,;\;P_2=\mathit i(P_1)\,;\; Q_1=\phi_1(P_1)\,;\; Q_2=\phi_2(P_2).$$
Indeed, let $\mathfrak p$ be a prime of $K$ of good reduction for $E$ and not over $3$ and suppose that $(m_1P_1+m_2P_2 \bmod \mathfrak p)$ has order coprime to $3$. It is clearly sufficient to show that both $(m_1P_1\bmod \mathfrak p)$ and $(m_2P_2\bmod \mathfrak p)$ have order coprime to $3$.
By multiplying $P_1$ and $P_2$ by an integer coprime to $3$, we may assume that $(P_1 \bmod \mathfrak p)=(R \bmod \mathfrak p)$ for a point $R$ in $E[3^{\infty}]$. Then we have $(m_1R+m_2i(R) \bmod \mathfrak p)=0$ and by the injectivity of the reduction modulo $\mathfrak p$ on $E[3^{\infty}]$ we deduce that $m_1R+m_2i(R)=0$. We have to show that $m_1R=0$. Let $3^h$ be the order of $R$. Then the annihilator of $R$ is an ideal of $\mathbb{Z}[i]$ containing $3^h$ but not $3^{h-1}$.
Since $3$ is prime in $\mathbb{Z}[i]$, the annihilator of $R$ is $(3^h)$.
Since $m_1 + m_2 i$ belongs to $(3^h)$, we can write $(m_1 + m_2 i) = 3^h(a_1 + a_2 i)$ for some integers $a_1$, $a_2$.
Therefore $m_1 R = 3^h a_1 R = 0$.
\end{exa}
 
We can also weaken condition (MSP) by imposing that $m_1=1$. Then one would like to prove that for every $i$ there exist $K$-endomorphisms $\phi_i$ and an integer $c$ such that $\phi_i(P_i)=cQ_i$. Without loss of generality it suffices to take $n=2$:
 
\begin{description}
\item[(WMSP)] Suppose that for all but finitely many primes $\mathfrak p$ of $K$ and for all positive integers $m$ the point $(Q_1+mQ_2 \bmod \mathfrak p)$ is zero whenever the point $(P_1+mP_2 \bmod \mathfrak p)$ is zero.
\end{description}

If $G$ is a simple abelian variety, under condition (WMPS) Bara\'nczuk showed that for $i=1,2$ there exist $K$-endomorphisms $\phi_i$ and an integer $c$ such that $\phi_i(P_i)=cQ_i$, see \cite[Theorem 8.1]{Baranczuk06}. The same proof holds for the multiplicative group hence for $1$-dimensional tori.
This result is in general false for a non-simple abelian variety or for a torus of dimension $>1$, as the following example shows.

\begin{exa}\label{nobar1}
Let $G$ be either an elliptic curve without complex multiplication or the multiplicative group defined over a number field $K$. Suppose that the rank of $G(K)$ is greater than $1$. Then let $(R,W)$ be a $K$-rational point on $G^2$ which is independent.
Consider the following points in $G^2(K)$:
$$P_1=Q_1=Q_2=(R,0)\,; P_2=(0,W).$$
These points satisfy condition (WMSP) but there do not exist a $K$-endomorphism $\phi$ of $G^2$ and a non-zero integer $c$ such that $\phi(P_2)=cQ_2$.
\end{exa}

\section*{Acknowledgements}

I thank Brian Conrad, Marc Hindry, Ren\'e Schoof for helpful discussions. I thank Jeroen Demeyer and Willem Jan Palenstijn for Example 16 and the case of $1$-dimensional tori. 


\vspace{0.4cm}

\noindent Math\'ematiques, \'Ecole Polytechnique F\'ed\'erale de Lausanne,\\ 1015 Lausanne, Switzerland
\newline 

\noindent \textit{e-mail}: antonella.perucca@epfl.ch

\end{document}